\documentclass{amsart}

\usepackage{amssymb,latexsym,amsmath,amsthm}

\newtheorem{thm}{Theorem}[section]

\theoremstyle{definition}

\theoremstyle{remark}

\numberwithin{equation}{section}

\begin{document}

\title{A RADEMACHER TYPE FORMULA FOR PARTITIONS AND OVERPARTITIONS}

\author{Andrew V. Sills}
\address{Department of Mathematical Sciences, 
Georgia Southern University, Statesboro, Georgia, 30460-8093 USA}
\email{ASills@GeorgiaSouthern.edu}



\date{\today}


\begin{abstract}
A Rademacher-type convergent series formula which generalizes the Hardy-Ramanujan-Rademacher formula for the number of partitions of $n$ and the Zuckerman formula for the
Fourier coefficients of $\vartheta_4(0 \mid \tau)^{-1}$ is presented. 
\end{abstract}

\maketitle

\newcommand{\F}{\mathcal{F}}
\newcommand{\I}{\mathcal{I}}

\section{Background}
\subsection{Partitions}
A \emph{partition} of an integer $n$ is a representation of $n$ as a sum of positive integers,
where the order of the summands (called \emph{parts}) is considered irrelevant.  It is customary to write the parts in nonincreasing order.
For example, there are three partitions of the integer $3$, namely $3$, $2+1$, and $1+1+1$.
Let $p(n)$ denote the number of partitions of $n$, with the convention that $p(0)=1$, and let $f(x)$ denote the generating
function of $p(n)$, i.e. let
\[ f(x) := \sum_{n=0}^\infty p(n) x^n .\]

Euler~\cite{le1748} was the first to systematically study partitions.  He showed that
\begin{equation}\label{egf}
f(x) = \prod_{m=1}^\infty \frac{1}{1-x^m} .
\end{equation}
Euler also showed that
\begin{equation} \label{epnt} \frac{1}{f(x)} = \sum_{n=-\infty}^\infty (-1)^n x^{n(3n-1)/2},
\end{equation}
and since the exponents appearing on the right side of~\eqref{epnt} are the pentagonal
numbers, Eq.~\eqref{epnt} is often called ``Euler's pentagonal number theorem."
  
   Although Euler's results can all be treated
from the point of view of formal power series, the series and infinite products above (and indeed all the series and infinite products mentioned in this paper) converge absolutely when $|x|<1$, which
is important for analytic study of these series and products.

Hardy and Ramanujan were the first to study $p(n)$ analytically and produced an
incredibly accurate asymptotic formula~\cite[p. 85, Eq. (1.74)]{hr1918}, namely
\begin{multline} \label{HRPofN}
 p(n) = \frac{1}{2\pi\sqrt{2}}
 \sum_{k=1}^{\lfloor \alpha\sqrt{n} \rfloor} \sqrt{k} 
 \underset{(h,k)=1}{\sum_{0\leqq h < k}} \omega(h,k) 
   e^{-2 \pi i h n/k} 
   \frac{d}{dn} \left(   \frac{ \exp\left( \frac{\pi}{k} \sqrt{\frac 23 \left( n-\frac{1}{24} \right)} \right) }
   { \sqrt{n-\frac{1}{24}}}\right) \\ + O(n^{-1/4}),
\end{multline}
where
\begin{equation*}
\omega(h,k) = \exp \left(  \pi i \sum_{r=1}^{k-1} \frac rk \left(  
\frac{hr}{k} - \lfloor \frac{hr}{k} \rfloor - \frac 12 \right) \right),
\end{equation*}
$\alpha$ is an arbitrary constant, 
and here and throughout $(h,k)$ is an abbreviation for $\gcd(h,k)$.

 Later Rademacher~\cite{har1937} improved upon~\eqref{HRPofN} by finding the following convergent series 
representation for $p(n)$:
\begin{equation} \label{RadPofN}
p(n) = \frac{1}{\pi\sqrt{2}} \sum_{k=1}^\infty \sqrt{k}
 \underset{(h,k)=1}
{ \sum_{0 \leqq h < k}}  \omega(h,k) e^{- 2\pi i nh/k}
\frac{d}{dn} \left( \frac{ \sinh \left( \frac{\pi}{k}  \sqrt{\frac23 \left( n-\frac{1}{24}\right)}\right)}{\sqrt{n-\frac{1}{24}}} \right).
\end{equation}  

Rademacher's method was used extensively by many practitioners,
including Grosswald~\cite{eg1958, eg1960}, Haberzetle~\cite{mh1941},
Hagis~\cite{ph1962, ph1963, ph1964a, ph1964b, ph1965a, ph1965b, ph1965c, ph1966, ph1971}, Hua~\cite{lkh1942}, Iseki~\cite{si1959, si1960, si1961}, Lehner~\cite{jl1941}, Livingood~\cite{jl1945}, 
Niven~\cite{in1940}, and Subramanyasastri~\cite{vvs1972}
to study various restricted partitions functions.

Recently, Bringmann and Ono~\cite{bo2009}
have given exact formulas for the coeffcients of all 
harmonic Maass forms of weight $\leqq \frac 12$.
The generating functions considered herein are 
weakly holomorphic modular forms of weight $-\frac 12$, 
and thus they are harmonic 
Maass forms of weight $\leqq \frac 12$. 
Accordingly, the results of this present paper could be derived from the general theorem 
in~\cite{bo2009}.  However, here we opt to derive the results via classical method of
Rademacher.
 
\subsection{Overpartitions} 
Overpartitions were introduced by S. Corteel and J. Lovejoy in~\cite{cl2004} and have been studied extensively by them and others including Bringmann, Chen, Fu, Goh, 
Hirschhorn, Hitczenko, Lascoux, Mahlburg, Robbins, R{\o}dseth, Sellers, Yee, and Zho
\cite{bl2007, czj2005, cgh2006, ch2004, cl2004, cly2004, cm2007, fl2005, hs2005, hs2006,
jl2003, jl2004a, jl2004b, jl2005a, jl2005b, jl2007, km2004, nr2003, rs2005}.

An~\emph{overpartition} of $n$ is a representation of $n$ as a sum of positive integers with summands in nonincreasing order, where the last occurrence of a given summand may or may not be overlined.  Thus the eight overpartitions of $3$ are $3$, $\bar{3}$, $2+1$, $\bar{2}+1$, $2+\bar{1}$, $\bar{2}+\bar{1}$, $1+1+1$, $1+1+\bar{1}$.  

Let $\bar{p}(n)$ denote the number of overpartitions of $n$
and let $\bar{f}(x)$ denote the generating function 
$\sum_{n=0}^\infty \bar{p}(n) x^n$ of $\bar{p}(n)$.  
Elementary techniques are sufficient to show that 
\[ \bar{f}(x) = \prod_{m=1}^\infty \frac{1+x^m}{1-x^m} = \frac{f(x)^2}{f(x^2)}. \]

 Note that 
 \[ \frac{1}{\bar{f}(x)} = \sum_{n=-\infty}^{\infty} (-1)^n x^{n^2} 
  \]
via an identity of Gauss~\cite[p. 23, Eq. (2.2.12)]{gea1976},
so that the reciprocal of the generating function for overpartitions is a series wherein a coefficient is nonzero if and only if the exponent of $x$ is a perfect square, just as the reciprocal of the generating function for partitions is a series wherein a coefficient is nonzero if and only if the exponent of $x$ is a pentagonal number.

Hardy and Ramanujan, writing more than 80 years before the coining of the term
``overpartition," stated~\cite[p. 109--110]{hr1918} that the function which we are calling
$\bar{p}(n)$ ``has no very simple arithmetical interpretation; but the series is none the less,
as the direct reciprocal of a simple $\vartheta$-funciton, of particular interest."  They went on
to state that 
\begin{equation} \label{HRPBarofN}
\bar{p}(n) = \frac{1}{4\pi} \frac{d}{dn}\left( \frac{ e^{\pi\sqrt{n}}}{\sqrt{n}} \right) +
   \frac{\sqrt{3}}{2\pi} \cos\left(  \frac 23 n \pi - \frac 16 \pi \right) \frac{d}{dn}  
    \left( e^{\pi\sqrt{n}/3} \right) + \cdots + O(n^{-1/4}).
\end{equation}  

In fact,~\eqref{HRPBarofN} was improved to the following Rademacher-type
convergent series by Zuckerman~\cite[p. 321, Eq. (8.53)]{hsz1939}:
\begin{equation}\label{overptnformula}
\bar{p}(n) = \frac{1}{2\pi} \underset{2\nmid k}{\sum_{k\geqq 1}} \sqrt{k}
 \underset{(h,k)=1}{\sum_{0\leqq h < k}} \frac{\omega(h,k)^2 }{ \omega(2h,k)  }
     e^{-2\pi i nh /k } \frac{d}{dn} \left(  \frac{\sinh \left( \frac{\pi \sqrt{n}}{k} \right) }{\sqrt{n} } \right).
 \end{equation}
A simplified proof of Eq.~\eqref{overptnformula} was given by L. Goldberg in this
Ph.D. thesis~\cite{lg1981}.

\subsection{Partitions where no odd part is repeated}
Let $pod(n)$ denote the number of partitions of $n$ where no odd part appears more than once.
Let $g(x)$ denote the generating function of $pod(n)$, so we have
\[ g(x) = \sum_{n=0}^\infty pod(n) x^n = \prod_{m=1}^\infty \frac{1+x^{2j-1}}{1-x^{2j}} 
= \frac{f(x)f(x^4)}{f(x^2)}. \]

Via another identity of Gauss~\cite[p. 23, Eq. (2.2.13)]{gea1976}, it turns out that
 \[ \frac{1}{g(x)} = \sum_{n=0}^\infty (-x)^{n(n+1)/2} = \sum_{n=-\infty}^\infty (-1)^n x^{2n^2-n}, \]
so in this case the reciprocal of the generating function under consideration has nonzero coefficients at the exponents which are triangular (or equivalently, hexagonal) numbers.

The analogous Rademacher-type formula for $pod(n)$ is as follows.
\begin{multline} \label{norepoddsformula}
 pod(n) = \frac{2}{\pi} \sum_{k\geqq 1}
 \sqrt{k \left( 1-(-1)^k + \lfloor \frac{(k,4)}{4} \rfloor  \right)} \\ \times
 \underset{(h,k)=1}{\sum_{0\leqq h < k}} \frac{\omega(h,k) \ \omega\left(  \frac{4h}{(k,4)}, \frac{k}{(k,4)} \right) }{ \omega\left( \frac{2h}{(k,2)} ,\frac{k}{(k,2)} \right)  }
     e^{-2\pi i nh /k } \frac{d}{dn} \left(  \frac{\sinh \left( \frac{\pi \sqrt{(k,4) (8n-1)}}{4k} \right) }{\sqrt{8n-1} } \right).
 \end{multline}

Eq.~\eqref{norepoddsformula} is the case $r=2$ of Theorem~\ref{MainFormula} to 
be presented in the next section.

\section{A common generalization}
Let us define 
\begin{equation} \label{frdef} f_r (x) :=  \frac{f(x) f(x^{2^r})}{f(x^2)}, \end{equation}
where $r$ is a nonnegative integer.
Thus,
\begin{align}
    f_0 (x) &=  \bar{f}(x) = \sum_{n=0}^\infty \bar{p}(n) x^n, \\
    f_1 (x) &=  f(x) = \sum_{n=0}^\infty p(n) x^n, \\
    f_2 (x) &= g(x) = \sum_{n=0}^\infty pod(n) x^n.
\end{align}

Let $p_r (n)$ denote the coefficient of $x^n$ in the expansion of $f_r(n)$, i.e.
\begin{equation*}
  f_r(x) = \sum_{n=0}^\infty p_r(n) x^n.
\end{equation*} 

Notice that $f_r(x) $ can be represented by several forms of equivalent infinite products,
each of which has a natural combinatorial interpretation:
\begin{align}
   f_r(x) & = \prod_{m=1}^\infty \frac{1+x^m}{1-x^{2^r m}}  \label{interp1}\\
      & = \prod_{m=1}^\infty \frac{1}{(1-x^{2m-1})(1-x^{2^r m})}\label{interp2}\\
      & = \prod_{m=1}^\infty \frac{1}{1-x^{2^{r-1}m} } \prod_{\lambda=1}^{2^{r-1}-1} (1+x^{2^{r-1}m+\lambda}).
      \label{interp3}
\end{align}

Thus, $p_r(n)$ equals each of the following:
\begin{itemize}
\item the number of overpartitions of $n$ where nonoverlined parts are multiples of $2^r$
(by~\eqref{interp1});
\item the number of partitions of $n$ where all parts are either odd or multiples of $2^r$
(by~\eqref{interp2}), provided $r\geqq 1$;
\item the number of partitions of $n$ where 
where nonmultiples of $2^{r-1}$ are distinct (by~\eqref{interp3}), provided $r\geqq 1$.
\end{itemize}

\begin{thm} \label{MainFormula}
For $r = 0, 1, 2, 3, 4$,  
\begin{multline*}
 p_r (n)  = \frac{2^{(r+1)/2} \sqrt{3}}{\pi}
 \underset{ (k, 2^{\max(r,1)})=1}{ \sum_{k\geqq 1} \sqrt{k}}
 \underset{(h,k)=1}{\sum_{0\leqq h < k}} e^{-2\pi i n h/k} \frac{\omega(h,k) \omega(2^r h, k)}
 {\omega(2h,k)} \\ \qquad\times \frac{d}{dn} \left\{ \frac{\sinh\left( 
  \frac{\pi \sqrt{(24n-2^r+1)(1+2^{r-1})}}{2^{r/2}\cdot 6k} \right)}{\sqrt{24n-2^r+1}} \right\}\\
 +  
 \frac{ \sqrt{3}}{\pi} \sum_{j=1+\lfloor \frac r2 \rfloor}^r 2^{(2-j+r)/2}
 \underset{ (k, 2^r)=2^j}{ \sum_{k\geqq 1} \sqrt{k}}
 \underset{(h,k)=1}{\sum_{0\leqq h < k}} e^{-2\pi i n h/k} \frac{\omega(h,k)
  \omega(2^{r-j} h, 2^{-j}k)}
 {\omega(h, \frac k2)} \\ \qquad \times \frac{d}{dn} \left\{ \frac{\sinh\left( 
  \frac{\pi \sqrt{(24n-2^r+1)(-1+2^{2j-r})}}{ 6k} \right)}{\sqrt{24n-2^r+1}} \right\}.
\end{multline*}
\end{thm}

 \section{A Proof of Theorem~\ref{MainFormula}} \label{OverPtnPf}
The method of proof is based on Rademacher's proof of~\eqref{RadPofN} 
in~\cite{har1943} with the necessary modifications. Additional details
of Rademacher's proof of~\eqref{RadPofN} are provided in \cite{har1955},
\cite[Ch. 14]{har1973} and \cite[Ch. 5]{tma1990}.

Of fundamental importance is the path of integration to be used.  
In~\cite{har1943}, Rademacher improved upon his original proof of~\eqref{RadPofN} given in~\cite{har1937},
by altering his path of integration from a carefully chosen circle to a more complicated path
based on Ford circles,
which in turn led to considerable simplifications later in the proof.  

\subsection{Farey fractions}

The sequence $\F_N$ of \emph{proper Farey fractions of order $N$} is the set of all $h/k$ with
$(h,k)=1$ and $0\leqq h/k <1$, arranged in increasing order.
Thus, e.g., $\F_4 = \left\{ \frac 01, \frac 14, \frac 13, \frac 12, \frac 23, \frac 34\right\}.$

For a given $N$, let $h_p$, $h_s$, $k_p$, and $k_s$ be such that $\frac{h_p}{k_p}$ is
the immediate predecessor of $\frac hk $ and $\frac{h_s}{k_s}$ is the immediate
successor of $\frac hk$ in $\F_N$.  It will be convenient to view each $\F_N$ cyclically, i.e.
to view
 $\frac 01$ as the immediate successor of
$\frac {N-1}{N}$.

\subsection{Ford circles and the Rademacher path}
Let $h$ and $k$ be integers with $(h,k) = 1$ and $0\leqq h < k$.  
The \emph{Ford circle}~\cite{lrf1938}
$C(h,k)$ is the circle in $\mathbb C$ of radius $\frac{1}{2k^2}$ centered at the point
$$ \frac{h}{k} + \frac{1}{2k^2} i.$$  

The \emph{upper arc $\gamma(h,k)$ of the Ford circle $C(h,k)$} is those points of $C(h,k)$
from the initial point 
\begin{equation} \label{AlphaI}
 \alpha_I(h,k):=  \frac hk - \frac{k_p}{k(k^2+k_p^2)} + \frac{1}{k^2+k_p^2} i
\end{equation}
to the terminal point 
\begin{equation}\label{AlphaT}
 \alpha_T(h,k):= \frac hk + \frac{k_s}{k(k^2+k_s^2)} + \frac{1}{k^2+k_s^2} i, 
 \end{equation}
traversed in the clockwise direction.

Note that we have
\[ \alpha_I(0,1)  = \alpha_T(N-1,N) . \]

Every Ford circle is in the upper half plane.
For $\frac{h_1}{k_1}, \frac{h_2}{k_2} \in \F_N$,  $C(h_1, k_1)$ and $C(h_2, k_2)$ are
either tangent or do not intersect.

The \emph{Rademacher path} $P(N)$ of order $N$ is the path in the upper half of the
$\tau$-plane from $i$ to $i+1$ consisting of 
\begin{equation} \label{RadPath}  \bigcup_{\frac hk \in \F_N}  \gamma(h,k) \end{equation}
traversed left to right and clockwise.  In particular, we consider the left half of the Ford
circle $C(0,1)$ and the corresponding upper arc $\gamma(0,1)$ to be translated
to the right by 1 unit.  This is legal given then periodicity of the function which is to 
be integrated over $P(N)$.


\subsection{Set up the integral}
Let $n$ and $r$ be fixed, with $n > (2^r - 1)/24$.

Since \[ f_r(x) = \sum_{n=0}^\infty p_r(n) x^n, \]
Cauchy's residue theorem implies that
\begin{equation} \label{Cauchy}
    p_r(n) = \frac{1}{2\pi i} \int_{\mathcal C} \frac{f_r(x)}{x^{n+1}}\ dx,
\end{equation}
where $\mathcal C$ is any 
simply closed contour enclosing the origin and inside the unit circle.  
   We introduce the change of variable
   \[ x = e^{2 \pi i \tau} \]
so that the unit disk $|x| \leqq 1$ in the $x$-plane maps to the infinitely tall, unit-wide strip in the
$\tau$-plane where $0 \leqq \Re\tau \leqq 1$ and $\Im\tau \geqq 0$.   
The contour $\mathcal C$ is then taken
to be the preimage of $P(N)$
under the map $x\mapsto e^{2\pi i \tau}$.

Better yet, let us replace $x$ with $e^{2 \pi i \tau}$ in~\eqref{Cauchy}
to express the integration in the $\tau$-plane:
\begin{align}\label{TauIntegral}
p_r(n) &= \int_{{P}(N)} {f_r}(e^{2\pi i \tau}) e^{-2\pi i n\tau} d\tau \nonumber\\
        &=  {\sum_{\frac hk \in \F_N}}
        \int_{\gamma(h,k)} {f_r}(e^{2\pi i \tau}) e^{-2\pi i n \tau} d\tau \nonumber\\
         &={\sum_{k=1}^{N}} \underset{(h,k)=1}{ \sum_{0 \leqq h < k}} \int_{\gamma(h,k)} 
       {f_r}(e^{2\pi i \tau}) e^{-2\pi i n \tau} d\tau
         \nonumber
\end{align}

\subsection{Another change of variable}
Next, we change variables again, taking
\begin{equation} \label{ZToTau} \tau = \frac {iz + h}{k}, \end{equation} so that
\begin{equation} \label{TauToZ} z = -ik \left(\tau - \frac hk\right). \end{equation}
Thus $C(h,k)$ (in the $\tau$-plane) maps to the clockwise-oriented 
circle $K^{(-)}_k$ (in the $z$-plane) centered
at $1/2k$ with radius $1/2k$.

So we now have 
\begin{multline}
p_r(n) 
=i{ \sum_{k=1}^{N}} k^{-1} \underset{(h,k)=1}
{ \sum_{0 \leqq h <  k}}    e^{-2\pi i n h/k} 
\underset{\mbox{\tiny arc of $K_k^{(-)} $} }
{\int_{z_I(h,k)}^{z_T(h,k)} }e^{2n\pi z/k} f_r( e^{2\pi i h/k - 2\pi z/k} ) \ dz
, \label{IntegralZ}
\end{multline}
where $z_I(h,k)$ (resp. $z_T(h,k)$) is the image of $\alpha_I(h,k)$ (see~\eqref{AlphaI})
(resp. $\alpha_T(h,k)$ [see~\eqref{AlphaT}]) under the transformation~\eqref{TauToZ}.

So the transformation~\eqref{TauToZ} maps the upper arc $\gamma(h,k)$ of $C(h,k)$
in the $\tau$-plane to the arc on $K^{(-)}_k$ which initiates at
\begin{equation} \label{ZI}
 z_I(h,k) = \frac{k}{k^2 + k_p^2} + \frac{ k_p}{k^2+ k_p^2} i
\end{equation} and terminates at
\begin{equation} \label{ZT}
  z_T(h,k) = \frac{k}{k^2+k_s^2} - \frac{ k_s}{k^2+k_s^2} i.
\end{equation}

\subsection{Exploiting a modular transformation}
From the theory of modular forms, we have the transformation formula~\cite[p. 93, Lemma 4.31]{hr1918}
\begin{multline} \label{FFunctionalEq}
 f\left(
    \exp \left( \frac {2\pi i h}{k} - \frac{2\pi z}{k} \right) \right)\\ = 
 \omega(h,k) \exp\left(\frac{\pi(z^{-1}-z)}{12k} \right)\sqrt{z} 
 f\left( \exp\left(2\pi i  \frac{ iz^{-1}+H}{k}\right) \right),
\end{multline}
where $\sqrt{z}$ is the principal branch, $(h,k)=1$, and $H$ is a solution to the congruence
\[ hH \equiv -1\pmod{k} . \]
From~\eqref{FFunctionalEq}, we deduce the analogous transformation for 
$f_r (x) $.





The transformation formula is a piecewise defined function with $r+1$ cases corresponding to $j=0,1,2,\dots, r$, where
$(k, 2^r) = 2^j$.
\begin{multline} \label{frFunctionalEq}
 f_r\left(
    \exp \left( \frac {2\pi i h}{k} - \frac{2\pi z}{k} \right) \right)\\ = 
\frac{ \omega(h,k) \omega(2^{r-j}h, 2^{-j} k ) }{ \omega\left( \frac{2h}{2-\delta_{j0}}, 
\frac{k}{2-\delta_{j0}} \right) }
\exp\left( \frac{\pi(2+2^{2j-r+1} - (2-\delta_{j0})^2)}{24kz} + \frac{\pi(1-2^r)z}{12k}    \right) \\
\times \sqrt{z \ 2^{r-j-1} (2-\delta_{j0})} 
\frac{ f\left( \exp\left(\frac{-2\pi}{kz} + \frac{2H_j\pi i}{k} \right)\right)  
   f\left( \exp\left(\frac{-2^{2j-r+1}\pi}{kz} + \frac{2^{2j-r+1}H_j\pi i}{k} \right)\right)           }
{ f\left( \exp\left(\frac{-\pi(2-\delta_{j0})^2}{kz} + \frac{ H_j\pi (2-\delta_{j0})^2 i}{k} \right)\right)    }
 ,
\end{multline}
where $H_j$ is divisible by $2^{r-j}$ and is a solution to the congruence 
$hH_j \equiv -1\pmod{k}$,
and 
\[ \delta_{j0} = \left\{ \begin{array}{ll} 1 & \mbox{if $j=0$} \\ 0 & \mbox{if $j\neq 0$}
\end{array} \right.
\] is the Kronecker $\delta$-function.

Notice that in particular, for $\lfloor \frac r2 \rfloor \leqq j \leqq r$,~\eqref{frFunctionalEq} simplifies to
\begin{multline} \label{frFunctionalEqSimp}
 f_r\left(
    \exp \left( \frac {2\pi i h}{k} - \frac{2\pi z}{k} \right) \right)\\ = 
\frac{ \omega(h,k) \omega(2^{r-j}h, 2^{-j} k ) }{ \omega\left( h, 
\frac{k}{2} \right) }
\exp\left( \frac{\pi}{12k} \left[  (2^{2j-r}-1)z^{-1} + (1-2^r)z\right]   \right) \\
\times \sqrt{z \ 2^{r-j} } 
 f_{2j-r}\left( \exp\left(\frac{-2\pi}{kz} + \frac{2H_j\pi i}{k} \right)\right) .
\end{multline}

Since the $r=0$ case was established
by Zuckerman, and the $r=1$ case by Rademacher, we will proceed with the 
assumption that $r>1$.

Apply~\eqref{frFunctionalEq} to~\eqref{IntegralZ} to obtain
\begin{multline}
p_r(n) 
=i \sum_{j=0}^{r} \underset{(k,2^r)=j}{ \sum_{k=1}^{N}} k^{-1} \underset{(h,k)=1}
{ \sum_{0 \leqq h <  k}}    e^{-2\pi i n h/k} 
\underset{\mbox{\tiny arc of $K_k^{(-)} $} }
{\int_{z_I(h,k)}^{z_T(h,k)} } 
\frac{ \omega(h,k) \omega(2^{r-j}h, 2^{-j} k ) }{ \omega\left( \frac{2h}{2-\delta_{j0}}, 
\frac{k}{2-\delta_{j0}} \right) } \\ \times
\exp\left( \frac{\pi(2+2^{2j-r+1} - (2-\delta_{j0})^2)}{24kz} + \frac{\pi(24n+1-2^r)z}{12k}    \right) \\
\times \sqrt{z \ 2^{r-j-1} (2-\delta_{j0})}  \\ \times
\frac{ f\left( \exp\left(\frac{-2\pi}{kz} + \frac{2H_j\pi i}{k} \right)\right)  
   f\left( \exp\left(\frac{-2^{2j-r+1}\pi}{kz} + \frac{2^{2j-r+1}H_j\pi i}{k} \right)\right)           }
{ f\left( \exp\left(\frac{-\pi(2-\delta_{j0})^2}{kz} + \frac{H_j\pi (2-\delta_{j0})^2 i}{k} \right)\right)    } \ dz
. 
\end{multline}

\subsection{Normalization} 
   Next, introduce a normalization $\zeta = zk$.  (This is not strictly necessary,
but it will allow us in the sequel to quote various useful results directly from the literature.)

\begin{multline}
p_r(n) 
=i \sum_{j=0}^{r} \underset{(k,2^r)=j}{ \sum_{k=1}^{N}} k^{-5/2} \underset{(h,k)=1}
{ \sum_{0 \leqq h <  k}}    e^{-2\pi i n h/k} 
\underset{\mbox{\tiny arc of $K^{(-)} $} }
{\int_{\zeta_I(h,k)}^{\zeta_T(h,k)} } 
\frac{ \omega(h,k) \omega(2^{r-j}h, 2^{-j} k ) }{ \omega\left( \frac{2h}{2-\delta_{j0}}, 
\frac{k}{2-\delta_{j0}} \right) } \\ \times
\exp\left( \frac{\pi(2+2^{2j-r+1} - (2-\delta_{j0})^2)}{24\zeta} + \frac{\pi(24n+1-2^r)\zeta}{12k^2}    \right) \\
\times \sqrt{\zeta \ 2^{r-j-1} (2-\delta_{j0})}  \\ \times
\frac{ f\left( \exp\left(\frac{-2\pi}{\zeta} + \frac{2H_j\pi i}{k} \right)\right)  
   f\left( \exp\left(\frac{-2^{2j-r+1}\pi}{\zeta} + \frac{2^{2j-r+1}H_j\pi i}{k} \right)\right)           }
{ f\left( \exp\left(\frac{-\pi(2-\delta_{j0})^2}{\zeta} + \frac{H_j\pi (2-\delta_{j0})^2 i}{k} \right)\right)    } 
\ d\zeta
, \label{IntegralZeta}
\end{multline}
where
\begin{equation} \label{ZetaI}
 \zeta_I(h,k) = \frac{k^2}{k^2 + k_p^2} + \frac{k k_p}{k^2+ k_p^2} i
\end{equation} and
\begin{equation} \label{ZetaT}
  \zeta_T(h,k) = \frac{k^2}{k^2+k_s^2} - \frac{k k_s}{k^2+k_s^2} i.
\end{equation}

Let us now rewrite~\eqref{IntegralZeta} as
\begin{multline} \label{RewrittenIntegralZeta}
p_r(n) 
=i \sum_{j=0}^{r} \underset{(k,2^{\max(r,1)})=2^j}{ \sum_{k=1}^{N}} k^{-5/2} \underset{(h,k)=1}
{ \sum_{0 \leqq h <  k}}    e^{-2\pi i n h/k} 
\frac{ \omega(h,k) \omega(2^{r-j}h, 2^{-j} k ) }{ \omega\left( \frac{2h}{2-\delta_{j0}}, 
\frac{k}{2-\delta_{j0}} \right) } \\ \times
 \left( \I_{j,1} + \I_{j,2} \right)
\end{multline}
where
\begin{multline} \I_{j,1}:=
\underset{\mathrm{arc}}{\int_{\zeta_I(h,k)}^{\zeta_T(h,k)} } 
\exp\left( \frac{\pi(2+2^{2j-r+1} - (2-\delta_{j0})^2)}{24\zeta} + \frac{\pi(24n+1-2^r)\zeta}{12k^2}    \right)  \\ \times \sqrt{\zeta \ 2^{r-j-1} (2-\delta_{j0})}
       \left\{ -1 +
       \frac{ f\left( \exp\left(\frac{-2\pi}{\zeta} + \frac{2H_j\pi i}{k} \right)\right)  
   f\left( \exp\left(\frac{-2^{2j-r+1}\pi}{\zeta} + \frac{2^{2j-r+1}H_j\pi i}{k} \right)\right)           }
{ f\left( \exp\left(\frac{-\pi(2-\delta_{j0})^2}{\zeta} + \frac{H_j\pi (2-\delta_{j0})^2 i}{k} \right)\right)    } 
       \right\} \ d\zeta,
\end{multline}
and
\begin{multline} \I_{j,2}:=
\underset{\mathrm{arc}}{\int_{\zeta_I(h,k)}^{\zeta_T(h,k)} }
\exp\left( \frac{\pi(2+2^{2j-r+1} - (2-\delta_{j0})^2)}{24\zeta} + \frac{\pi(24n+1-2^r)\zeta}{12k^2}    
\right) \\ \times \sqrt{\zeta \ 2^{r-j-1} (2-\delta_{j0})} \ d\zeta.
\end{multline}

\subsection{Estimation}  
It will turn out that as $N\to\infty$, only $\I_{j,2}$ for $j=0$ and $\lfloor r/2 \rfloor < j \leqq r$ ultimately make a contribution (provided $r<5$).
Note that all the integrations in the $\zeta$-plane occur on arcs and chords of the circle $K$
of radius $\frac 12$ centered at the point $\frac 12$.   So, inside and on $K$,  $0<\Re\zeta\leqq 1$
and $\Re \frac 1\zeta \geqq 1$.

\subsubsection{Estimation of $\I_{j,2}$ for $1\leqq j \leqq \lfloor r/2 \rfloor$}
The regularity of the integrand allows us to alter the path of integration from the arc connecting
$\zeta_I(h,k)$ and $\zeta_T(h,k)$ to the directed segment.
By \cite[p. 104, Thm. 5.9]{tma1990}, the length of the path of integration does not exceed $2\sqrt{2} k/N$, and on the segment connecting $\zeta_I(h,k)$ to $\zeta_T(h,k)$, $|\zeta|< \sqrt{2} k /N$.
Thus, the absolute value of the integrand,
   \begin{align*}
 & \phantom{=} \left|  \exp\left( \frac{\pi(2^{2j-r} - 1)}{12\zeta} + \frac{\pi(24n+1-2^r)\zeta}{12k^2}    
\right) \sqrt{\zeta \ 2^{r-j}} \right| \\
   &= |\zeta|^{1/2} 2^{(r-j)/2} \exp\left(\frac{(24n+1-2^r)\pi \Re\zeta}{12k^2}  \right)
   \exp\left( \frac{\pi(2^{2j-r}-1 )}{12} \Re\frac{1}{\zeta} \right)\\
   & \leqq |\zeta|^{1/2} 2^{r/2} \exp(2\pi n) .
   \end{align*}
Thus,  for $1\leqq j \leqq \lfloor r/2 \rfloor$,
\[ | \I_{j,2} | \leqq \frac{2\sqrt{2} k}{N} \left( \frac{\sqrt{2} k }{N} \right)^{1/2} 2^{r/2} e^ {2\pi n} \leqq
C'_j k^{3/2} N^{-3/2}, \]
for a constant $C'_j$ (recalling that $n$ and $r$ are fixed).

\subsubsection{Estimation of $\I_{j,1}$ for $1\leqq j \leqq \lfloor r/2 \rfloor$}
 We have the absolute value of the integrand:
\begin{align*}
  & \left| \sqrt{\zeta \ 2^{r-j}} \exp\left( \frac{\pi(2^{2j-r} - 1)}{12\zeta} + \frac{\pi(24n+1-2^r)\zeta}{12k^2}    \right)  \right| \\ & \qquad \times \left|
  -1 +
       \frac{ f\left( \exp\left(\frac{-2\pi}{\zeta} + \frac{2H_j\pi i}{k} \right)\right)  
   f\left( \exp\left(\frac{-2^{2j-r+1}\pi}{\zeta} + \frac{2^{2j-r+1}H_j\pi i}{k} \right)\right)           }
{ f\left( \exp\left(\frac{-4\pi}{\zeta} + \frac{4 H_j\pi i}{k} \right)\right)    } 
       \right| \\
=& \left| \sqrt{\zeta \ 2^{r-j}} \exp\left( \frac{\pi(2^{2j-r} - 1)}{12\zeta} + \frac{\pi(24n+1-2^r)\zeta}{12k^2}    \right)  \right| \\ & \qquad \times \left|
  -1 + f_{2j-r} \left( \exp\left(  \frac{-2\pi}{\zeta} + \frac{2H_j\pi i}{k} \right) \right)
       \right| \\       
= &|\zeta|^{1/2} 2^{(r-j)/2} \exp\left( \frac{ (24n+1-2^r) \pi  \Re\zeta}{12k^2} \right)
\exp\left(  \frac{(2^{2j-r}-1)\pi}{12} \Re \frac 1\zeta \right)\\ & \qquad \times  
\left| \sum_{m=1}^\infty p_{2j-r} (m) 
\exp\left( \frac{-2\pi m}{\zeta} + \frac{2H_j\pi i m}{k}  \right)\right| \\
\leqq &|\zeta|^{1/2} 2^{(r-j)/2} \exp\left( \frac{ (24n+1-2^r) \pi }{12} \right)
\sum_{m=1}^\infty p_{2j-r} (m) 
\exp\left(  -\frac{\pi}{12} (24m - 2^{2j-r}+ 1) \right) \\
\leqq &|\zeta|^{1/2} 2^{r/2} e^{2\pi n}
\sum_{m=1}^\infty p_{0} (m) e^{-2\pi m} \\
= & c_j |\zeta|^{1/2},
\end{align*} for a constant $c_j$.
So, for $1\leqq j \leqq \lfloor r/2 \rfloor$,
 \[ |\I_{j,1}| \leqq \frac{2\sqrt{2} k}{N} \left( \frac{\sqrt{2} k }{N} \right)^{1/2} c_j 
 < C_j k^{3/2} N^{-3/2} \]
for a constant $C_j$.

\subsubsection{Estimation of $I_{0,1}$}
Let $p^*_r(x)$ be defined by
\[  \sum_{n=0}^\infty p^*_r(n) x^n =  \frac{ f(x^{2^r}) f(x)}{ f(x^{2^{r-1}}) }. \]

Again, the regularity of the integrand allows us to alter the path of integration from the arc connecting
$\zeta_I(h,k)$ and $\zeta_T(h,k)$ to the directed segment.

With this in mind, we estimate the absolute value of the integrand:
\begin{align*}  
&\phantom{=} \left | \exp\left( \frac{\pi(1+2^{1-r}) }{24\zeta} + \frac{\pi(24n+1-2^r)\zeta}{12k^2}    \right)  \sqrt{\zeta \ 2^{r-1} } \right| \\ &\qquad\times
       \left| -1 +
       \frac{ f\left( \exp\left(\frac{-2\pi}{\zeta} + \frac{2H_0\pi i}{k} \right)\right)  
   f\left( \exp\left(\frac{-2^{1-r}\pi}{\zeta} + \frac{2^{1-r}  H_0\pi i}{k} \right)\right)           }
{ f\left( \exp\left(\frac{-\pi}{\zeta} + \frac{H_0\pi i}{k} \right)\right)    } 
       \right| \\
&= \left | \exp\left( \frac{\pi(1+2^{1-r}) }{24\zeta} + \frac{\pi(24n+1-2^r)\zeta}{12k^2}    \right)  \sqrt{\zeta \ 2^{r-1} } \right| \\ &\qquad\times
       \left| \sum_{m=1}^\infty p^*(m) \exp \left(  \frac{-2^{1-r} \pi m}{\zeta} + 
       \frac{2^{1-r}H_0 \pi i m}{k} \right)
       \right| \\
&= \exp\left( \frac{\pi(1+2^{1-r}) }{24} \Re\frac{1}{\zeta} \right) 
\exp\left( \frac{\pi(24n+1-2^r)\Re\zeta}{12k^2}    \right)  
|\zeta|^{1/2} 2^{(r-1)/2}   \\ &\qquad\times
       \left| \sum_{m=1}^\infty p^*(m) \exp \left(  \frac{-2^{1-r} \pi m}{\zeta} \right)
     \exp\left( \frac{2^{1-r}H_0 \pi i m}{k} \right)
       \right|   \\
&\leqq e^{2\pi n}
|\zeta|^{1/2} 2^{(r-1)/2}   \\ &\qquad\times
       \sum_{m=1}^\infty  | p^*(m) | 
        \exp\left( \frac{\pi(1+2^{1-r}) }{24} \Re\frac{1}{\zeta}-2^{1-r} \pi m\Re\frac{1}{\zeta} \right)\\
&= e^{2\pi n}
|\zeta|^{1/2} 2^{(r-1)/2}   
       \sum_{m=1}^\infty  | p^*(m) | 
        \exp\left( -\frac{\pi }{24\cdot 2^{r-1}} \Re\frac{1}{\zeta} \left( 24m - 1 - 2^{r-1}  \right) \right)   \\
 &\leqq e^{2\pi n}
|\zeta|^{1/2} 2^{(r-1)/2}   
       \sum_{m=1}^\infty  | p^*(m) | 
        \exp\left( -\frac{\pi }{24} \left( 24m - 1 - 2^{r-1}  \right) \right)     \\
 &< e^{2\pi n}
|\zeta|^{1/2} 2^{(r-1)/2}   
       \sum_{m=1}^\infty  | p^*(24m-1-2^{r-1}) | 
       y^{24m - 1 - 2^{r-1}}  \quad\mbox{(where $y=e^{-\pi/24}$)}      \\
    & =c_0 |\zeta|^{1/2},       
    \end{align*}
for a constant $c_0$.
So, \[ |\I_{0,1}| \leqq \frac{2\sqrt{2} k}{N} \left( \frac{\sqrt{2} k }{N} \right)^{1/2} c_0 
< C_0 k^{3/2} N^{-3/2} \]
for a constant $C_0$.

\subsubsection{Estimation of $\I_{j,1}$ for $1+\lfloor r/2 \rfloor \leqq j \leqq r \leqq 4$}
\label{EstIj1-2}

Again, the regularity of the integrand allows us to alter the path of integration from the arc connecting
$\zeta_I(h,k)$ and $\zeta_T(h,k)$ to the directed segment.

With this in mind,
\begin{align}
&\phantom{=}\left| \exp\left( \frac{\pi(2^{2j-r} -1) }{12\zeta} + \frac{\pi(24n+1-2^r)\zeta}{12k^2}    \right) \sqrt{\zeta \ 2^{r-j} } \right| \notag \\ 
& \qquad\times 
       \left| -1 +
        f_{2j-r}\left( \exp\left(\frac{-2\pi}{\zeta} + \frac{2H_j\pi i}{k} \right)\right)  
          \right|  \notag\\
&= \left| \exp\left( \frac{\pi(2^{2j-r} -1) }{12\zeta} \right)
\exp\left( \frac{\pi(24n+1-2^r)\zeta}{12k^2}    \right) \sqrt{\zeta \ 2^{r-j} } \right| \notag\\ 
& \qquad\times 
       \left| \sum_{m=1}^\infty p_{2j-r}(m) \exp\left( \frac{-2\pi m}{\zeta} \right)
       \exp\left(\frac{2H_j\pi i m}{k} \right)       \right| \notag  \\
&=  \exp\left( \frac{\pi(2^{2j-r} -1) }{12}\Re\frac{1}{\zeta} \right)
\exp\left( \frac{\pi(24n+1-2^r)\Re\zeta}{12k^2}    \right) |\zeta|^{1/2}  2^{(r-j)/2} \notag\\ 
& \qquad\times 
       \left| \sum_{m=1}^\infty p_{2j-r}(m) \exp\left( {-2\pi m} \Re\frac{1}{\zeta} \right)
       \exp\left(\frac{2H_j\pi i m}{k} \right)       \right|  \notag \\
&\leqq  
e^{2\pi n} |\zeta|^{1/2}  2^{r/2}
  \sum_{m=1}^\infty  p_{2j-r}(m) 
  \exp\left( 
    -\frac{\pi}{12}\Re\frac{1}{\zeta} (24m- 2^{2j-r} +1 )  \right)  \label{later}\\     
&\leqq  
e^{2\pi n} |\zeta|^{1/2}  2^{r/2}
  \sum_{m=1}^\infty  p_{0}(m) 
  \exp\left( 
    -\frac{\pi}{12}\Re\frac{1}{\zeta} (24m- 2^{r} +1 )  \right)\notag\\                            
&=  c_j |\zeta|^{1/2} \notag
\end{align} for a constant $c_j$.
So, \[ |\I_{j,1}| \leqq \frac{2\sqrt{2} k}{N} \left( \frac{\sqrt{2} k }{N} \right)^{1/2} c_j 
< C_j k^{3/2} N^{-3/2} \]
for a constant $C_j$, when $1+\lfloor r/2 \rfloor\leqq j \leqq r$.

\subsubsection{Combining the estimates}
\begin{align*}
&\left| 
i \underset{(k,2^r)=1}{\sum_{k=1}^N} k^{-5/2} \underset{(h,k)=1}{ \sum_{0\leqq h<k}}
  e^{-2\pi i n h/k} \frac{\omega(h,k) \omega(2^r h,k)}{\omega(2h,k)}
 \I_{0,1} \right.  \\
&\left. + i \sum_{j=1}^{\lfloor r/2 \rfloor}
 \underset{(k,2^r)=2^j}{\sum_{k=1}^N} k^{-5/2} \underset{(h,k)=1}{ \sum_{0\leqq h<k}}
  e^{-2\pi i n h/k} \frac{\omega(h,k) \omega(2^{r-j}h, 2^{-j}k)}
  {\omega(h,k/2)} \left( \I_{j,1} + \I_{j,2} \right) \right. \\
&\left. + i \sum_{j=1+\lfloor r/2 \rfloor}^{r}
 \underset{(k,2^r)=2^j}{\sum_{k=1}^N} k^{-5/2} \underset{(h,k)=1}{ \sum_{0\leqq h<k}}
  e^{-2\pi i n h/k} \frac{\omega(h,k) \omega(2^{r-j}h, 2^{-j}k)}
  {\omega(h,k/2)}   \I_{j,1} \right| \\  
  &< \sum_{j=0}^r \sum_{k=1}^N \sum_{h=0}^{k-1} C_j k^{-1}  N^{-3/2} 
  + \sum_{j=1}^{\lfloor r/2 \rfloor}\sum_{k=1}^N \sum_{h=0}^{k-1}  C'_j k^{-1} N^{-3/2}\\
 & \leqq C'' N^{-3/2} \sum_{k=1}^N 1,\qquad \mbox{(where $C'' = \sum_{j=0}^r C_j
 + \sum_{j=1}^{\lfloor r/2 \rfloor} C'_j$)}  \\
  &= O(N^{-1/2}).
\end{align*}
Thus, we may revise~\eqref{IntegralZeta} to
\begin{multline}
p_r(n) =
i \sum_{k=1}^N k^{-5/2} \underset{(h,k)=1}{\sum_{0\leqq h<k} } e^{-2\pi i n h/k}
\frac{\omega(h,k) \omega(2^{r-j}h, 2^{-j} k)}{\omega(2h,k)} \I_{0,2} \\
+ i \sum_{j=1+\lfloor \frac r2 \rfloor}^r 
\sum_{k=1}^N k^{-5/2} \underset{(h,k)=1}{\sum_{0\leqq h<k} } e^{-2\pi i n h/k}
\frac{\omega(h,k) \omega(2^{r-j}h, 2^{-j} k)}{\omega(h,k/2)} \I_{j,2}
 + O(N^{-1/2}).
 \label{IntegralZetaO}
\end{multline}

\subsection{Evaluation of $\I_{j,2}$ for $1+\lfloor \frac r2 \rfloor \leqq j \leqq r $}
  Write $\I_{j,2}$ as
\begin{multline} \I_{j,2} = \int_{K^{(-)}} 
\exp\left( \frac{\pi(2^{2j-r} -1)}{12\zeta} + \frac{\pi(24n+1-2^r)\zeta}{12k^2}    
\right)  \sqrt{\zeta \ 2^{r-j} } \ d\zeta \\
- \I_{j,3}- \I_{j,4}, \label{FullCircle}\end{multline} 
where
\[ \I_{j,3} := \int_{0}^{\zeta_I(h,k)}, 
\quad \I_{j,4} := \int_{\zeta_T(h,k)}^{0},
\] 
and $\I_{j,3}$ and $\I_{j,4}$ have the same integrand as~\eqref{FullCircle}.
\subsubsection{Estimation of $\I_{j,3}$ and $\I_{j,4}$}
We note that the length of the arc of integration in $\I_{j,3}$ 
is less than $\frac{\pi k}{\sqrt{2}N}$,
and on this arc $|\zeta| < \sqrt{2} k / N$.~\cite[p. 272]{har1973}.
Also, $\Re \frac 1\zeta = 1$ on $K$~\cite[p. 271, Eq. (120.2)]{har1973}.
Further, $0<\Re \zeta < 2k^2/N^2$~\cite[p. 271, Eq. (119.6)]{har1973}.
The absolute value of the integrand is thus
\begin{align*}
&|2^{r-j} \zeta |^{1/2}  \exp\left( \frac{ (24n+1-2^r) \pi \Re\zeta}{12k^2} 
+ \frac{(2^{2j-r}-1)\pi}{12}\Re \frac 1\zeta\right)\\
<& 2^{(r-j)/2} 2^{1/4} k^{1/2} N^{-1/2} 
\exp\left( \frac{(24n+1-2^r)\pi}{6 N^2} + \frac{(2^{2j-r}-1)\pi}{12} \right)
\end{align*} so that
\begin{align*}  |\I_{j,2}| 
&< \pi k 2^{-1/2} N^{-1} 2^{(r-j)/2} 2^{1/4} k^{1/2} N^{-1/2}
\exp\left( \frac{(24n+1-2^r)\pi}{6 N^2} + \frac{(2^{2j-r}-1)\pi}{12} \right) \\
&=  \pi k^{3/2} N^{-3/2} 
 2^{(2r-2j-1)/4} 
\exp\left( \frac{(24n+1-2^r)\pi}{6 N^2}  + \frac{(2^{2j-r}-1)\pi}{12} \right) \\ 
 &= O\left(k^{3/2} N^{-3/2} \exp\left( \frac{(24n+1-2^r)\pi}{6 N^2}  
 \right) \right)  . \end{align*}
By the same reasoning, $|\I_{j,3}| =
O\left(k^{3/2} N^{-3/2} \exp\left( \frac{(24n+1-2^r)\pi}{6 N^2}  
 \right) \right)  . $

 We may therefore revise~\eqref{IntegralZetaO} to
\begin{multline}
p_r(n) =
i \sum_{k=1}^N k^{-5/2} \underset{(h,k)=1}{\sum_{0\leqq h<k} } e^{-2\pi i n h/k}
\frac{\omega(h,k) \omega(2^{r-j}h, 2^{-j} k)}{\omega(2h,k)} \\ \times
\int_{K^{(-)}}  
\sqrt{\zeta 2^{r-j-1}}
\exp\left\{  \frac{\pi(2^{1-r}+1)}{24\zeta} + \frac{\pi(24n+1-2^r)\zeta}{12k^2}
\right\} \ d\zeta\\
+ i \sum_{j=1+\lfloor \frac r2 \rfloor}^r 
\sum_{k=1}^N k^{-5/2} \underset{(h,k)=1}{\sum_{0\leqq h<k} } e^{-2\pi i n h/k}
\frac{\omega(h,k) \omega(2^{r-j}h, 2^{-j} k)}{\omega(h,k/2)} \\ \times
\int_{K^{(-)}}  
\sqrt{\zeta 2^{r-j}}
\exp\left\{  \frac{\pi(2^{2j-r}-1)}{12\zeta} + \frac{\pi(24n+1-2^r)\zeta}{12k^2}
\right\} \ d\zeta
 + O(N^{-1/2})
 \label{IntegralZetaO2}
\end{multline}
and upon letting $N$ tend to infinity, obtain
\begin{multline}
p_r(n) =
i \sum_{k=1}^\infty k^{-5/2} \underset{(h,k)=1}{\sum_{0\leqq h<k} } e^{-2\pi i n h/k}
\frac{\omega(h,k) \omega(2^{r-j}h, 2^{-j} k)}{\omega(2h,k)} \\ \times
\int_{K^{(-)}}  
\sqrt{\zeta 2^{r-j-1}}
\exp\left\{  \frac{\pi(2^{1-r}+1)}{24\zeta} + \frac{\pi(24n+1-2^r)\zeta}{12k^2}
\right\} \ d\zeta\\
+ i \sum_{j=1+\lfloor \frac r2 \rfloor}^r 
\sum_{k=1}^\infty k^{-5/2} \underset{(h,k)=1}{\sum_{0\leqq h<k} } e^{-2\pi i n h/k}
\frac{\omega(h,k) \omega(2^{r-j}h, 2^{-j} k)}{\omega(h,k/2)} \\ \times
\int_{K^{(-)}}  
\sqrt{\zeta 2^{r-j}}
\exp\left\{  \frac{\pi(2^{2j-r}-1)}{12\zeta} + \frac{\pi(24n+1-2^r)\zeta}{12k^2}
\right\} \ d\zeta
 \label{IntegralZetaNoO}
\end{multline}

\subsection{The final form}
We may now introduce the change of variable 
\[ \zeta = \frac{\pi\left( (2+2^{2j-r+1}-(2-\delta_{j0})^2) \right)}{24t}, \]
(where the first summation in~\eqref{IntegralZetaNoO} is the $j=0$ term separated
out for clarity),
which allows the integral to be evaluated
in terms of $I_{3/2}$, the Bessel function of the first kind of order $3/2$ with purely imaginary argument~\cite[p. 372, \S17.7]{ww1927} when we bear in mind that a ``bent" path of 
integration is allowable according to the remark preceding Eq. (8) on p. 177 of~\cite{gnw1944}.
See also~\cite[p. 109]{tma1990}.
The final form of the formula is then obtained by using the fact that Bessel functions of half-odd integer order can be expressed in terms of elementary functions.

We therefore have
\begin{multline*} 
p_r(n) = \frac{\pi}{(24n-2^r+1)^{3/4}} \left\{  
  \frac{(1+2^{r-1})^{3/4}}{2^{(r-2)/4}} \underset{(k,2^{\max(r,1)})=1)}{\sum_{k\geqq 1}}
  k^{-1} \right. \\ \times \underset{(h,k)=1}{\sum_{0\leqq h<k}} e^{-2\pi i n h/k} 
  \frac{\omega(h,k) \omega(2^r h,k)}{\omega(2h,k)} I_{3/2} 
  \left( \frac{\pi \sqrt{(24n-2^r+1)(1+2^{r-1})}}{2^{r/2}\cdot6k} \right) \\
  + \sum_{j=1+\lfloor \frac r2 \rfloor}^r (2^{2j-r}-1 )^{3/4} 2^{(2-j+r)/2}
  \underset{(k,2^r)=2^j}{\sum_{k\geqq 1}}
  k^{-1}\\ \left.
   \times \underset{(h,k)=1}{\sum_{0\leqq h<k}} e^{-2\pi i n h/k} 
  \frac{\omega(h,k) \omega(2^{r-j} h, 2^{-j}k)}{\omega(h,k/2)} I_{3/2} 
  \left( \frac{\pi \sqrt{(24n-2^r+1)(1+2^{2j-r})}}{6k} \right) \right\}
  \end{multline*}
which, after application of the formula~\cite[p. 110]{tma1990}
\begin{equation}
I_{3/2} (z) = \sqrt{ \frac{2z}{\pi}} \frac{d}{dz} \left( \frac{\sinh z}{z} \right),
\end{equation}
 is equivalent to Theorem~\ref{MainFormula}. \hfill $\Box$
  
\section*{Acknowledgment}
The author thanks the anonymous referee for bringing the work of Zuckerman~\cite{hsz1939}
and Goldberg~\cite{lg1981} to his attention.  This, in turn, led the author to seek the more
general result presented here in this final version of the paper.

\end{document}